\documentstyle[emlines2,bezier,12pt]{article}
\begin{document}
\title{On some cases of integrability of a general Riccati equation}
\author{N.M.~Kovalevskaya\\
{\small
{\it Novgorod State University}}}

\date{}

\maketitle

\begin{center}
{\bf Abstract}
\end{center}

{\footnotesize

A general Riccati equation is integrated in quadratures in case
one of its coefficients is an arbitrary function and two others
are expressed through it.

}
\vskip 1.5cm

It is known that a general Riccati equation
\begin{equation}
\label{sys1}
y'=P(t)y^{2}+Q(t)y+F(t)
\end{equation}
can be integrated in quadratures only in special cases. These ones
are given, for instance, in [1,2].

We will show that it is possible to find the solution of the
equation~(\ref{sys1}) in case one of the coefficients P,Q or F is
an arbitrary function and two others can be expressed through it.

Changing $\displaystyle y=-\frac{1}{P}\frac{u'(t)}
{u(t)}$ we can
write the equation~(\ref{sys1}) in the form of an equation of the
second order in regards to a new variable $u$:
\begin{equation}
\label{sys2}
u''-(\frac{P'}{P}+Q)u'+PFu=0
\end{equation}

Now we consider equation~(\ref{sys2}) as a system of two equations
of the first order:
\begin{eqnarray}
\label{sys3} \left \{
\begin{array}{rcl}
u'&=&z  \\
z'&=&-PFu+(\displaystyle\frac{P'}{P}+Q)z
\end{array}
\right.
\end{eqnarray}

that can be rewritten as a matrix equation $\dot X=A(t)X$ (here
$X=\left(
\begin{array}{c}
u \\
v
\end{array}
\right) )$ with the matrix of coefficients:
\begin{equation}
\label{sys4}
 A(t)=
\left(
\begin{array}{cc}
 0 & 1  \\
 -PF & \frac{P'}{P}+Q
\end{array}
\right).
\end{equation}

It is also known [3], that if a matrix $A(t)$ is of the form

\begin{equation}
\label{sys5} \left(
\begin{array}{cc}
a_{11}(t)  &  a_{12}(t)  \\
c_{1}a_{12}(t)  &  a_{11}(t)+c_{2}a_{12}(t)

\end{array}
\right)
\end{equation}
where $c_{1}$ and $c_{2}$ are arbitrary constants then $A(t)$
satisfies the functional commutativity condition
\begin{equation}
\label{6} A(t')A(t'')-A(t'')A(t')= 0, \quad
   \forall t',t''
\in{\bf R}
\end{equation}
and a fundamental matrix of solutions of a corresponding system is
the matrix exponential $X(t)=\exp\left(\int\limits_0^t
A(\tau)d\tau\right)$ [3].

Now we prove that under a special choice of the elements of the
matrix~(\ref{sys4}) it is the functionally commutative one so the
system~(\ref{sys3}) is integrable and the equation~(\ref{sys2})
has an explicit solution.

Considering $A(t)$ of the form~(\ref{sys5}) we have $a_{11}(t)=
0$, $a_{12}(t)=1$. Due to relations for elements $a_{21}(t)$ and
$a_{22}(t)$ we get a system for determination of coefficients of
the Riccati equation:
\begin{eqnarray}
\label{sys7} \left \{
\begin{array}{rcl}
c_{1}&=&-PF  \\
c_{2}&=&\displaystyle\frac{P'}{P}+Q.
\end{array}
\right.
\end{eqnarray}

We will show that solving of~(\ref{sys1}) is possible when one of
the coefficients F,P or Q is an arbitrary one and two others are
expressed through it.

Let $F(t)$ be an arbitrary function and $c_{1}\neq 0$ in
~(\ref{sys7}), then we have
$$
P(t)=\frac{-c_{1}}{F(t)}; \quad
Q(t)=c_{2}-\frac{P'(t)}{P(t)}=c_{2}+\frac{F'(t)}{F(t)}.
$$

We substitute the expressions for $P(t)$ and $Q(t)$
into~(\ref{sys4}) and find
\begin{equation}
\label{sys8} A(t)= \left(
\begin{array}{cc}
0 & 1  \\
c_{1} & c_{2} \\
\end{array}
\right).
\end{equation}

So the matrix of the system~(\ref{sys3}) turns out to be a matrix
with constant coefficients and a corresponding Riccati
equation~(\ref{sys1}) can be written as follows:

\begin{equation}
\label{sys9}
y'=-\frac{c_{1}}{F(t)}y^2+(c_{2}+\frac{F'(t)}{F(t)})y+F(t).
\end{equation}

If we take $P(t)$ as an arbitrary function in~(\ref{sys7}) then
functions $F(t)$ and $Q(t)$ can be immediately expressed through
it and they have the form
$$
F(t)=-\frac{c_{1}}{P(t)}; \quad  Q(t)=c_{2}-\frac{P'(t)}{P(t)}.
$$

With these expressions for the functions we have the same
form~(\ref{sys8}) of $A(t)$ and the initial Riccati equation gives
over into
\begin{equation}
\label{sys10}
y'=P(t)y^2+(c_{2}-\frac{P'(t)}{P(t)})y-\frac{c_{1}}{P(t)}.
\end{equation}

If we express $P(t)$ and $F(t)$ through $Q(t)$ then we obtain
$\displaystyle \frac{P'(t)}{P(t)}=c_{2}-Q(t)$ and
$P(t)=C\exp^{\displaystyle c_{2}t}\exp^{\displaystyle \int
Q(t)dt}.$

For $F(t)$ we have the following equality:
$$
F(t)=-\frac{c_{1}}{P(t)}=-\frac{c_{1}}{C} \exp^{-c_{2}t}\exp^{\int
Q(t)dt}.
$$

In this case $A(t)$ also has the form~(\ref{sys8}) and the Riccati
equation can be written as
\begin{equation}
\label{sys11} y'=C\exp^{\displaystyle
c_{2}t}\exp^{\displaystyle{-\int Q(t)dt}}
y^{2}+Q(t)y-\frac{c_{1}}{C}\exp^{\displaystyle
{-c_{2}t}}\exp^{\displaystyle \int Q(t)dt}.
\end{equation}

Thus, for functionally commutative matrix of the
system~(\ref{sys3}) which has the form~(\ref{sys5}) there are
three possible variants of the Riccati equation depending on our
choice of arbitrary functions $P(t), F(t)$ or $Q(t)$. All three
equations are integrable in quadratures.

For instance, if $F(t)$ is an arbitrary function then the solution
of~(\ref{sys9}) is of the form:
$$
y(t)=\frac{((2c_{1}+Cc_{2})\sinh(\frac{1}{2}\sqrt{c_{2}^{2}+4c_{1}}t)+C\sqrt{c_{2}^{2}+4c_{1}}\cosh(\frac{1}{2}\sqrt{c_{2}^{2}+4c_{1}}t))F(t)}
{((2C-c_{2})\sinh(\frac{1}{2}\sqrt{c_{2}^{2}+4c_{1}}t)+\sqrt{c_{2}^{2}+4c_{1}}\cosh(\frac{1}{2}\sqrt{c_{2}^{2}+4c_{1}}t))c_{1}}
$$

where $C$ is an arbitrary constant which is determined from the
initial conditions.

Letting $P(t)$ be an arbitrary function we obtain the following
form of solution of~(\ref{sys10}):
$$
y(t)=\frac{(2c_{1}+Cc_{2})\sinh(\frac{1}{2}\sqrt{c_{2}^{2}+4c_{1}}t)+C\sqrt{c_{2}^{2}+4c_{1}}\cosh(\frac{1}{2}\sqrt{c_{2}^{2}+4c_{1}}t)}
{((c_{2}-2C)\sinh(\frac{1}{2}\sqrt{c_{2}^{2}+4c_{1}}t)-\sqrt{c_{2}^{2}+4c_{1}}\cosh(\frac{1}{2}\sqrt{c_{2}^{2}+4c_{1}}t))P(t)}.
$$

If choose $Q(t)$ as an arbitrary function then we get the solution
of ~(\ref{sys11}) in the form:
\begin{eqnarray*}
y(t) &=& \exp^{\displaystyle -c_{2}t}\exp^{\displaystyle {\int
Q(t)dt}} \times \\
     & \times& \frac{((2c_{1}+Cc_{2})\sinh(\frac{1}{2}\sqrt{c_{2}^{2}+4c_{1}}t)+C\sqrt{c_{2}^{2}+4c_{1}}\cosh(\frac{1}{2}\sqrt{c_{2}^{2}+4c_{1}}t))}
{(c_{2}-2C)\sinh(\frac{1}{2}\sqrt{c_{2}^{2}+4c_{1}}t)-\sqrt{c_{2}^{2}+4c_{1}}\cosh(\frac{1}{2}\sqrt{c_{2}^{2}+4c_{1}}t)A}.
\end{eqnarray*}

where $A$ can be assumed to equal $1$ without loss of generality.

These found cases of the solution in quadratures of the general
Riccati equation are different from the known ones [2] and they
can be useful when solving some applied problems.

\par

\bigskip

\begin{thebibliography}{3}
\bibitem{1} Matveev N. M. {\it Methods of Ordinary Differential Equations Integration.} M.,1963;\\
\bibitem{2} Kamke E. {\it Ordinary Differential Equations Handbook.} M., 1961; \\
\bibitem{3} Kovalevskaya N.M., Kovalevsky M.M. {\it The conditions of
integrability in quadratures the second order linear systems of
differential equations.}, dep. RISTI, N 2450-B00, 2000.
\end{thebibliography}
\end{document}